\documentclass[reqno,11pt]{article}

\usepackage[utf8]{inputenc}
\usepackage[T1]{fontenc}

\usepackage{amsmath}   
\usepackage{amssymb}  
\usepackage{amsthm}    
\usepackage{mathrsfs}  
\usepackage{mathtools} 
\usepackage{bm}

\usepackage{enumitem}  
\usepackage{graphicx} 
\usepackage{color}     
\usepackage[dvipsnames]{xcolor}
\usepackage{cite} 
\usepackage{authblk} 
\usepackage{hyperref}
\hypersetup{
	colorlinks=true,   
	linkcolor=blue,     
	citecolor=teal,    
	urlcolor=cyan,      
	pdftitle={}, 
	pdfauthor={}
}
\usepackage{cleveref}   

\usepackage[a4paper]{geometry}

\theoremstyle{plain} 
\newtheorem{theorem}{Theorem}
\newtheorem{lemma}{Lemma}
\newtheorem{corollary}{Corollary}
\newtheorem{proposition}{Proposition}

\theoremstyle{definition}
\newtheorem{definition}{Definition}
\newtheorem{example}{Example}

\Crefname{condition}{Condition}{Conditions}

\theoremstyle{remark}
\newtheorem{remark}{Remark}

\begin{document}

\title{Grid complexity for two-dimensional words} 
\author[a]{Cheng Zeng \thanks{czeng@sdtbu.edu.cn}}
\author[b]{Qingcheng Zeng \thanks{Corresponding author: qczeng@buaa.edu.cn}}
\author[b]{Yumei Xue \thanks{yxue@buaa.edu.cn}}

\affil[a]{{\small{ School of Mathematics and Information Science, Shandong Technology and Business University, 264003, P. R. China}}}
\affil[b]{{\small{School of Mathematical Science, Beihang University, Beijing  100191, P. R. China}}}

\date{}
\maketitle
\vspace{-2em}

\begin{abstract}
	For a two-dimensional word \(\alpha\), let \(G_\alpha(m,n)\) be the
	largest number of patterns occurring on translates of a Cartesian
	product of \(m\) horizontal and \(n\) vertical coordinates. We compare
	this grid complexity with rectangular block complexity and maximal
	pattern complexity. We prove that \(G_\alpha(m,n)=2mn\) for binary
	interval codings of circle rotations with two irrational coordinate
	increments and for the Kamae--Xue word. For binary simple Toeplitz
	words, the same formula holds for all \(m,n\) if and only if neither coordinate
	direction is periodic. We also relate grid complexity to sequence
	entropy. For binary words, positive grid entropy is equivalent to
	\(G_\alpha(m,n)=2^{mn}\) for every \(m,n\). The Pascal word attains
	this maximum although its rectangular complexity $R_P(n,n)$ is \(O(n^2)\).
	
\end{abstract}
\section{Introduction}

Let \(\mathcal A\) be a finite alphabet. An element \(\alpha=\alpha_0\alpha_1\alpha_2\ldots\in\mathcal A^\mathbb N\) is called  a \emph{one-sided word} over \(\mathcal A\) if all letters in \(\mathcal A\) appear in \(\alpha\). For a finite window \(E\) of size \(k\), we mean a strictly increase integer sequence of length \(k\) with
\[
e_0=0<e_1<e_2<\cdots<e_{k-1}.
\]
We also refer to a window of size $k$ as a window with $k$ sites.
We say \(E\) is a block window of size $k$ if \(E=\{0,1,\ldots,k-1\}\). Denote by \(p_\alpha(E)=\#\{\alpha_{n+e_0}\alpha_{n+e_1}\cdots\alpha_{n+e_{k-1}}:\ n\in \mathbb N\}\) the
number of distinct patterns occurring on translates of \(E\). We denote the block complexity \(b_\alpha(k):= p_\alpha(E)\) for a block window of size $k$.
\emph{Maximal pattern complexity} \(p_\alpha^*(k)\) is the maximum of \(p_\alpha(E)\) over all windows of size \(k\). 

For one-sided infinite words, a well-known result by Morse and Hedlund \cite{mh} says that
\begin{theorem}
The following statements are equivalent:
\begin{enumerate}
\renewcommand{\labelenumi}{\textup{(\roman{enumi})}}
\item \(\alpha\) is eventually periodic;
\item \(b_\alpha(k)\) is bounded in \(k\);
\item \(b_\alpha(k)<k+1\) for some $k\in \mathbb Z_+$.
\end{enumerate}
\end{theorem}

Kamae and Zamboni \cite{kz} proved a fundamental result for maximal pattern complexity.
\begin{theorem}
	The following statements are equivalent:
	\begin{enumerate}
		\renewcommand{\labelenumi}{\textup{(\roman{enumi})}}
		\item \(\alpha\) is eventually periodic;
		\item \(p_\alpha^*(k)\) is bounded in \(k\);
		\item \(p_\alpha^*(k)<2k\) for some $k\in \mathbb Z_+$.
	\end{enumerate}
\end{theorem}
For recent results on one-dimensional maximal pattern complexity, see \cite{schlortt,lps}. Kamae, Rao and Xue \cite{krx} studied the corresponding periodicity problem for two-dimensional words under a recurrence hypothesis. 

For a word \(\alpha\in A^{\mathbb Z^2}\), we consider windows of the form \(S\times T\), where \(S,T\) are windows in \(\mathbb Z\). We define the \emph{grid complexity} of \(\alpha\) by
\[
G_\alpha(m,n)=
\max_{\substack{S,T\subset\mathbb Z\\ |S|=m,\ |T|=n}}
p_\alpha(S\times T).
\]
Thus the horizontal and vertical coordinates may be chosen
independently, and the chosen grid is translated as a whole.
Writing \(R_\alpha(m,n)\) for the number of distinct blocks on a
consecutive \(m\times n\) rectangle \(\{0,1,\ldots,m-1\}\times\{0,1,\ldots,n-1\}\), we have
\[
R_\alpha(m,n)\le G_\alpha(m,n)\le p_\alpha^*(mn).
\]
The second inequality raises a question: when can a Cartesian-product window attain the maximal pattern complexity? Based on the work of Kamae et al. \cite{kx, krx} and Sander and Tijdeman \cite{ST2002}, we establish several criteria for grid complexity and periodicity.

Kamae and Xue \cite{kx} constructed a two-dimensional binary word with \(p_\alpha^*(k)=2k\) for every \(k\). Qu, Rao, Wen and
Xue \cite{qrwx} obtained the same formula for binary simple Toeplitz words. For this latter class, we prove that
\begin{align*}
&G_\alpha(m,n)=2mn \quad\text{for all }m,n\ge1
\\
&\quad\Longleftrightarrow
\alpha\text{ has neither a horizontal nor a vertical period}.
\end{align*}

We also consider the full  grid complexity for a binary word, namely \(2^{mn}\). One example attaining this value is the Pascal word obtained from Pascal's triangle modulo \(2\).
Allouche and Berth\'e \cite{AB1997} studied the rectangular block complexity of Pascal's triangle modulo an integer. For the binary Pascal word, we exhibit an \(m\times n\) grid whose horizontal
translates realize every binary pattern on that grid.

The paper is organized as follows. Section \ref{sec:pre} establishes the basic properties of grid complexity. Section \ref{sec:bounded} gives periodicity criteria for grids of fixed width and for grids of arbitrary size. Section \ref{sec:simple} proves the criterion for simple Toeplitz words and discusses the Kamae--Xue word. Section \ref{sec:irr} proves \(G_\alpha(m,n)=2mn\) for binary rotation words when both coordinate increments are irrational. Finally, Section \ref{sec:full} gives two constructions with full grid complexity, one based on binary digits and the other on Pascal's triangle.

\section{Preliminaries}\label{sec:pre}
\subsection{Notation and basic facts}
Let $\alpha\in \mathcal A^{\mathbb Z^2}$ be a two dimensional word, where the letter set $\mathcal A$ satisfies $|\mathcal A|=q\ge2$. We also write \(A=\mathcal A\). Let \(\alpha(u)\) denote the letter of the word \(\alpha\) at position \(u=(i,j)\in\mathbb Z^2\). We adopt the convention for the shift \((\sigma^z\alpha)(u)=\alpha(u+z)\). Let $X_\alpha=\overline{\{\sigma^z\alpha:z\in\mathbb Z^2\}}$ be the subshift generated by $\alpha$, and we call $X_\alpha$ the orbit closure of $\alpha$.

For a finite window $E\subset\mathbb Z^2$, we define the factor set and its cardinality by
\[
\mathcal F_\alpha(E)=\{(\alpha(z+e))_{e\in E}:z\in\mathbb Z^2\},\quad p_\alpha(E)=|\mathcal F_\alpha(E)|.
\]
Here patterns retain their coordinate labels. Since finite cylinders are open, we have $\mathcal F_\alpha(E)=\{x|_E:x\in X_\alpha\}$.

From the definition of $p_\alpha(E)$, we can easily check the following.
\begin{proposition}\label{prop:mono}
For two finite windows $E,F$, we have 
\begin{enumerate}
		\renewcommand{\labelenumi}{\textup{(\roman{enumi})}}
		\item $p_\alpha(E)\leq p_\alpha(F)$ if $E\subset F$;
		
		\item $p_\alpha(E\cup F)\leq p_\alpha(E)p_\alpha(F)$.
\end{enumerate}
\end{proposition}

\begin{definition}[Three complexities]
Let $I_m=\{0,1,\ldots,m-1\}$ be the block window of size $m$. For positive $m,n,N$, set
\begin{align}
	R_\alpha(m,n)&=p_\alpha(I_m\times I_n),\\
	G_\alpha(m,n)&=\max_{\substack{S,T\subset\mathbb Z\\|S|=m,\ |T|=n}}p_\alpha(S\times T),\\
	p_\alpha^*(N)&=\max_{\substack{E\subset\mathbb Z^2\\|E|=N}}p_\alpha(E).
\end{align}
We call $R_\alpha(m,n),	G_\alpha(m,n), p_\alpha^*(N)$ rectangular, grid, and maximal pattern complexity, respectively.
\end{definition}

\begin{proposition}
All three complexities are nondecreasing in their size variables, and
\begin{align}
		R_\alpha(m,n)&\le G_\alpha(m,n)\le p_\alpha^*(mn)\le q^{mn},\label{chain}\\
		p_\alpha^*(N)&\le G_\alpha(N,N)\le p_\alpha^*(N^2).\label{completion}
\end{align}
If the two coordinate projections of $E$ have sizes $m,n$, then $p_\alpha(E)\leq G_\alpha(m,n)$.
\end{proposition}
\begin{proof}
For windows of fixed size $m\times n$, the first inequality follows since rectangular windows are special grid windows, which in turn are special general windows. 
	
The inequality $p^*(N)\leq G(N,N)$ holds since for any general window $E$ of size $N$, we can choose a grid window $S\times T$ of size $N\times N$ such that $E\subset S\times T$. The inequality $G(N,N)\leq p^*(N^2)$ follows from \eqref{chain}.

If the two coordinate projections of $E$ have sizes $m,n$, then there exists a window $S\times T$ such that $E\subset S\times T$, which implies that $p_\alpha(E)\leq p_\alpha(S\times T)\leq G_\alpha(m,n)$.
\end{proof}
\begin{proposition}\label{prop:submulti}
The complexities $R_\alpha$, $G_\alpha$ and $p_\alpha^*$ satisfy the sub-multiplicative property, that is, for positive integers $m,n,m_1,m_2,n_1,n_2$, 
\begin{align*}
R_\alpha(m_1+m_2,n)&\le R_\alpha(m_1,n)R_\alpha(m_2,n),\\
R_\alpha(m,n_1+n_2)&\le R_\alpha(m,n_1)R_\alpha(m,n_2),\\
G_\alpha(m_1+m_2,n)&\le G_\alpha(m_1,n)G_\alpha(m_2,n),\\
G_\alpha(m,n_1+n_2)&\le G_\alpha(m,n_1)G_\alpha(m,n_2),\\
p_\alpha^*(k+\ell)&\le p_\alpha^*(k)p_\alpha^*(\ell).
\end{align*}
\end{proposition}
\begin{proof}
Proposition \ref{prop:submulti} follows directly from Proposition \ref{prop:mono}.
\end{proof}
\begin{proposition}\label{prop:differ}
If $\alpha,\beta$ differ at at most $d$ sites, then
\[|p_\alpha(E)-p_\beta(E)|\le d|E|.\]
Consequently, $|R_\alpha(m,n)-R_\beta(m,n)|\leq dmn$, $|G_\alpha(m,n)-G_\beta(m,n)|\leq dmn$ and $|p_\alpha^*(N)-p_\beta^*(N)|\leq dN$. 
\end{proposition}
\begin{proof}
Let $D=\{u\in\mathbb Z^2:\alpha(u)\neq \beta(u)\}$, so that $|D|\leq d$. Write $P_x(z,E)=(x(z+e))_{e\in E}$ for the pattern read at $z+E$.

The two patterns $P_\alpha(z,E)$ and $P_\beta(z,E)$ differ precisely when $z+e=u$ for some $e\in E$ and $u\in D$. Thus the set of affected translations is 
\[
B_E=D-E=\{u-e: u\in D, e\in E\},
\]
which means that each pair $(u,e)$ determines exactly one such translation $z=u-e$. Different pairs may yield the same translation, so $|B_E|\leq |D||E|\leq d|E|$. 

Consider a pattern $w\in\mathcal F_\alpha(E)\setminus\mathcal F_\beta(E)$ and choose a
translation $z$ with $P_\alpha(z,E)=w$. If $z\notin B_E$, then
$P_\beta(z,E)=P_\alpha(z,E)=w$, contrary to the choice of $w$.
Hence every such $w$ occurs among the $\alpha$-patterns read at translations in
$B_E$. Each translation yields just one $\alpha$-pattern, giving
\[
|\mathcal F_\alpha(E)\setminus\mathcal F_\beta(E)|\le |B_E|.
\]
Now count the patterns common to both families, then the ones present only in
$\alpha$. This gives $p_\alpha(E)\le p_\beta(E)+|B_E|$.
Interchanging $\alpha$ and $\beta$ proves
$|p_\alpha(E)-p_\beta(E)|\le d|E|$.

These three complexity statements are immediate corollaries.
\end{proof}
\begin{remark}
The upper bound $d|E|$ in Proposition \ref{prop:differ} is attainable. Let $\alpha=0^{\mathbb Z^2}$ and $\beta=\bm{1}_{u=(0,0)}$. Then $R_\alpha(m,n)=G_\alpha(m,n)=p_\alpha^*(mn)=1$ and $R_\beta(m,n)=G_\beta(m,n)=p_\beta^*(mn)=mn+1$. Here $d=1$, and each pair of corresponding complexities differs by exactly $mn$.
\end{remark}
A word is \emph{strongly recurrent} if every centered square has two exact
return vectors whose angle stays a fixed positive distance from both $0$
and $\pi$, uniformly over square sizes. \emph{Uniform recurrence} means that
every occurring finite pattern has a relatively dense set of occurrences.
Uniform recurrence implies strong recurrence. By Kamae, Rao and Xue
\cite[Theorem 1.5]{krx}, a strongly recurrent word that is not doubly
periodic satisfies $p_\alpha^*(k)\ge2k$ for every $k\ge1$.

\subsection{Grid entropy and sequence entropy}

We define the \emph{grid entropy} of \(\alpha\) by
\[
h_{\mathrm{grid}}(\alpha)
=\lim_{m,n\to\infty}\frac{\log G_\alpha(m,n)}{mn},
\]
where both \(m\) and \(n\) tend to infinity. To see that the limit
exists, fix \(r,s\ge1\) and put \(u=\lfloor m/r\rfloor\) and
\(v=\lfloor n/s\rfloor\). Any \(m\times n\) grid can be partitioned
into \(uv\) grids of size \(r\times s\) and \(mn-uvrs\) remaining
sites. Hence
\[
G_\alpha(m,n)
\le G_\alpha(r,s)^{uv}\,q^{mn-uvrs}.
\]
Taking logarithms, dividing by \(mn\), and then taking the infimum
over \(r,s\) gives
\[
h_{\mathrm{grid}}(\alpha)
=\inf_{r,s\ge1}\frac{\log G_\alpha(r,s)}{rs}
=\lim_{k\to\infty}\frac{\log G_\alpha(k,k)}{k^2}.
\]
In short, the existence of this limit follows from the sub-multiplicativity of grid complexity.

Recall that \(X_\alpha\) is the shift orbit closure of \(\alpha\). Let
\(\mathcal P_0\) be the partition of \(X_\alpha\) according to the symbol at the
origin. For a sequence \(\mathbf g=(g_j)_{j\ge1}\) of distinct sites
in \(\mathbb Z^2\), the topological sequence entropy of this partition
along \(\mathbf g\) is given by
\[
h_{\mathbf g}(X_\alpha,\mathcal P_0)
=\limsup_{N\to\infty}
\frac{1}{N}\log p_\alpha(\{g_1,\ldots,g_N\}).
\]
Indeed, each pattern on these sites determines one nonempty atom of the corresponding joined partition; see~\cite{goodman} for
topological sequence entropy.
\begin{remark}
Huang and Ye~\cite{hy} relate maximal pattern entropy to sequence entropy over unrestricted sampling sequences. The proposition below uses the partition \(\mathcal P_0\) and sequences whose first \(k^2\) sites form a grid for every \(k\).
\end{remark}
\begin{proposition}
	The grid entropy \(h_{\mathrm{grid}}(\alpha)\) is the maximum of
	\(h_{\mathbf g}(X_\alpha,\mathcal P_0)\) over all sequences
	\(\mathbf g\) for which
	\[
	\{g_1,\ldots,g_{k^2}\}=S_k\times T_k \qquad(k\ge1),
	\]
	where \(S_k=\{s_1,\ldots,s_k\}\) and \(T_k=\{t_1,\ldots,t_k\}\)
	for some strictly increasing integer sequences \((s_i)\) and \((t_i)\).
\end{proposition}

\begin{proof}
	For such a sequence, the pattern count at \(N=k^2\) is at most
	\(G_\alpha(k,k)\). Since prefix pattern counts are nondecreasing
	and \((k+1)^2/k^2\to1\), taking the limsup along square indices yields
	the same sequence entropy as taking it along all indices. Thus
	\(h_{\mathbf g}(X_\alpha,\mathcal P_0)\le h_{\mathrm{grid}}(\alpha)\).
	
	To attain equality, choose integers \(L_j\) such that
	\((L_1+\cdots+L_{j-1})/L_j\to0\). For each \(j\), take a grid
	\(A_j\times B_j\) attaining \(G_\alpha(L_j,L_j)\). Translate its
	two coordinate sets so that each lies beyond all previously chosen
	coordinates on its axis. This preserves its pattern count.
	Concatenate the ordered sets \(A_j\) and \(B_j\) on their respective
	axes to obtain \((s_i)\) and \((t_i)\). At \(K_j=L_1+\cdots+L_j\),
	the grid \(S_{K_j}\times T_{K_j}\) contains the \(j\)-th maximizing
	grid. It therefore gives
	\[
	\frac{\log p_\alpha(S_{K_j}\times T_{K_j})}{K_j^2}
	\ge
	\left(\frac{L_j}{K_j}\right)^2
	\frac{\log G_\alpha(L_j,L_j)}{L_j^2}
	\longrightarrow h_{\mathrm{grid}}(\alpha).
	\]
	Enumerating the sites added at each successive grid gives the
	required sequence \(\mathbf g\).
\end{proof}

For a binary word, the supports of its patterns on \(S\times T\) are exactly the traces of the translates of \(\{z:\alpha(z)=1\}\) on that grid. Thus the binary assertion below also follows from the box shatter bound of Chernikov, Palac\'in and Takeuchi~\cite[Proposition~3.9]{cpt}: if no \(k\times k\) grid realizes every pattern for some \(k\), then \(\log G_\alpha(N,N)=o(N^2)\).

\begin{theorem}\label{thm:grid-entropy-spectrum}
	Let $2\le d\le q$. The following conditions are equivalent.
	\begin{enumerate}
		\renewcommand{\labelenumi}{\textup{(\roman{enumi})}}
		\item $h_{\rm grid}(\alpha)>\log(d-1)$.
		\item There is a set $B\subseteq A$ of cardinality $d$ such that,
		for every $m,n\ge1$, there exist $S,T\subset\mathbb Z$ with
		$|S|=m$, $|T|=n$, and
		\[
		B^{S\times T}\subseteq\mathcal F_\alpha(S\times T).
		\]
	\end{enumerate}
	Consequently,
	\[
	h_{\rm grid}(\alpha)\in\{\log 1,\log 2,\ldots,\log q\}.
	\]
	In particular, for a binary word,
	\[
	h_{\rm grid}(\alpha)>0
	\quad\Longleftrightarrow\quad
	G_\alpha(m,n)=2^{mn}\quad\text{for every }m,n\ge1.
	\]
\end{theorem}

\begin{proof}
	Condition (\romannumeral2) gives $G_\alpha(m,n)\ge d^{mn}$ for all $m,n$,
	and hence $h_{\rm grid}(\alpha)\ge\log d$. We prove the converse
	by contraposition.
	
	We first recall the counting bound of Daniely, Schapira and Shahaf \cite[Theorem~1.5]{dss}. For a nonempty family
	$\mathcal H\subseteq A^E$, let $v$ be the largest cardinality of
	a set $D\subseteq E$ for which one can choose $d$-element sets
	$B_e\subseteq A$, $e\in D$, satisfying
	\[
	\prod_{e\in D} B_e\subseteq\mathcal H|_D.
	\]
	The sets $B_e$ may differ from site to site. With $M=|E|$ and
	$b=\binom qd$, we have
\begin{equation}\label{eq:count}
	|\mathcal H|
\le \sum_{j=0}^{v}\binom Mj b^j(d-1)^{M-j}.
\end{equation}
To see this, delete one coordinate and let $\mathcal H_a$ be the
	family of prefixes admitting the letter $a$ there. A prefix with
	$t\ge1$ extensions satisfies $t\le d-1+\binom td$. Thus
	\[
	|\mathcal H|
	\le (d-1)\Bigl|\bigcup_{a\in A}\mathcal H_a\Bigr|
	+\sum_{\substack{C\subseteq A\\ |C|=d}}
	\Bigl|\bigcap_{a\in C}\mathcal H_a\Bigr|.
	\]
	The union has parameter at most $v$. Each nonempty intersection has
	parameter at most $v-1$, since the deleted coordinate can be restored
	with all letters of $C$. If $v=0$, all these intersections are empty.
	Induction on $M$, starting at $M=0$, now gives  \eqref{eq:count} by Pascal's
	identity.
	
	Suppose that (\romannumeral2) fails. Each $d$-element set $C\subseteq A$ then
	fails at some grid size. Since there are only finitely many such
	$C$, there is a common $r\ge2$ such that no $r\times r$ grid
	realizes all patterns over any one of them. Indeed, a larger full
	grid restricts to a full grid of every smaller size.
	
	Fix a $k\times k$ grid $S\times T$ and apply \eqref{eq:count} to
	$\mathcal H=\mathcal F_\alpha(S\times T)$. Regard a set $D$
	attaining $v$ as a bipartite graph with vertex classes $S$ and $T$.
	Color each edge $e$ by its chosen set $B_e$. A copy of $K_{r,r}$
	in color $C$ would give an $r\times r$ subgrid on which every
	assignment over $C$ occurs. Hence each color class is $K_{r,r}$-free.
	The K\H{o}v\'ari--S\'os--Tur\'an theorem \cite{kst} gives
	\[
	v=O_{q,r}(k^{2-1/r}),
	\]
	where the implied constant in \(O_{q,r}(\cdot)\) depends only on $q$ and $r$.
	Substituting $M=k^2$ into \eqref{eq:count}, uniformly over the chosen grid,
	yields
	\[
	\log G_\alpha(k,k)
	\le k^2\log(d-1)+O_{q,r}(k^{2-1/r}\log k).
	\]
	Dividing this inequality by $k^2$ and letting $k\to\infty$ implies $h_{\rm grid}(\alpha)\le\log(d-1)$, proving the equivalence of (\romannumeral1) and (\romannumeral2). 
	
	Since (\romannumeral2) implies
	$h_{\rm grid}(\alpha)\ge\log d$, no entropy value lies strictly
	between $\log(d-1)$ and $\log d$. The binary assertion follows
	by taking $d=q=2$.
\end{proof}

Every value in the theorem is attained. Fix $1\le d\le q$ and
choose a $d$-element set $B\subseteq A$. Place every finite square
block over $B$ in infinitely many pairwise disjoint squares and
fill the remaining sites with one letter of $B$. Then every
consecutive $m\times n$ rectangle realizes all $d^{mn}$ patterns.
If all $q$ letters are required to occur, insert each of the $q-d$
remaining letters at one site. Every old pattern still occurs,
while at most $(q-d)mn$ translates of a fixed $m\times n$ grid
meet an altered site. Thus
\[
d^{mn}\le G_\alpha(m,n)\le d^{mn}+(q-d)mn,
\]
and $h_{\rm grid}(\alpha)=\log d$.
\section{Bounded complexity and periodicity}\label{sec:bounded}
Set $L(\alpha)=\{z:\sigma^z\alpha=\alpha, z\in \mathbb Z^2\}$. The word $\alpha$ has a nonzero period when $\operatorname{rank}L(\alpha)\ge1$. It is \emph{doubly periodic} when $\operatorname{rank}L(\alpha)=2$, equivalently when $L(\alpha)$ has finite index in $\mathbb Z^2$.

\begin{theorem}
The following are equivalent: 
\begin{enumerate}
	\renewcommand{\labelenumi}{\textup{(\roman{enumi})}}
	\item $\alpha$ is doubly periodic;
	\item $X_\alpha$ is finite;
	\item $\sup_Np^*(N)$ is finite;
	\item $\sup_{m,n}G_\alpha(m,n)$ is finite;
	\item $\sup_{m,n}R_\alpha(m,n)$ is finite.
\end{enumerate}
If $X_\alpha$ is infinite, then
\[p_\alpha^*(N)\geq N+1\]
holds for all $N\geq 1$. The latter bound is sharp without recurrence assumptions.
\end{theorem}
\begin{proof}
The orbit is in bijection with $\mathbb Z^2/L(\alpha)$, and a finite orbit is closed. A finite orbit closure bounds every complexity. 

Conversely, if $X_\alpha$ contained $C+1$ distinct points, finitely many coordinates would distinguish them. A rectangle containing those coordinates would have at least $C+1$ patterns. Thus bounded block complexity forces a finite orbit closure $X_\alpha$. The same holds for the other two complexity quantities. 

Thus, we prove the equivalence of the five statements.

For the last claim, start with the empty window and one atom $X_\alpha$. After finitely many cylinder refinements, some atom is still infinite. Two points in it differ at a new coordinate, whose addition increases the number of atoms by at least one. After $N$ steps there are at least $N+1$ patterns.
\end{proof}

\begin{theorem}
For any fixed $w\ge1$,
\[
\sup_n R_\alpha(w,n)<\infty\ \Longleftrightarrow\
\sup_n G_\alpha(w,n)<\infty\ \Longleftrightarrow\
(0,b)\in L(\alpha)\text{ for some }b\ge1.
\]
The transposed assertion holds for a horizontal period.
\end{theorem}
\begin{proof}
If $b$ is a common vertical period, then $G_\alpha(w,n)\le q^{wb}<\infty$, and $R_\alpha(w,n)\le G_\alpha(w,n)$.

Conversely, suppose $R_\alpha(w,n)\leq C$ for all $n\geq 1$. Hence $R_\alpha(1,n)\leq C$. We now show that this bounds the number of distinct bi-infinite columns, even when the vertical origin may be chosen independently in each column. Let
\[
Y=\{y_{i,j}:i,j\in\mathbb Z\},\quad y_{i,j}(t)=\alpha(i,j+t), \quad t\in\mathbb Z.
\]
If $Y$ contained $C+1$ distinct sequences $y_1,\ldots,y_{C+1}$,
then for each pair $r<s$ we could choose a coordinate $t_{rs}$ at which
$y_r$ and $y_s$ differ. There are finitely many pairs. Take a single
integer interval $[L,U]$ containing all the chosen coordinates.
The restrictions of the $C+1$ sequences to $[L,U]$ are pairwise distinct.
Each restriction is a vertical word of length $U-L+1$ occurring in $\alpha$:
for $y_{i,j}$ its first site is $(i,j+L)$. This gives
$R_\alpha(1,U-L+1)\ge C+1$, a contradiction. Thus $|Y|\le C$.

It remains to obtain one period that works in every column. The vertical
shift $V$, defined by $(Vy)(t)=y(t+1)$, sends $y_{i,j}$ to $y_{i,j+1}$.
Its inverse sends $y_{i,j}$ to $y_{i,j-1}$, so $V$ permutes the finite set $Y$.
Let $b$ be the least common multiple of the lengths of its cycles. Then $V^b$
fixes every member of $Y$. In particular, $V^by_{i,0}=y_{i,0}$ for every $i$.
Evaluating at $j$ gives $\alpha(i,j+b)=\alpha(i,j)$ for all $i,j\in\mathbb Z$.
This proves the vertical assertion. Interchanging the two coordinates proves
the horizontal assertion.
\end{proof}
\begin{example}\label{ex:contex}
Let $u$ be a two-sided Sturmian word and put $\alpha(i,j)=u(i)$. 
This word is uniformly recurrent and its period group $L(\alpha)$ is $\{0\}\times\mathbb Z$. 
Its one-dimensional block and maximal pattern counts mean that
\[
R_\alpha(m,n)=m+1,\quad G_\alpha(m,n)=2m,\quad p_\alpha^*(k)=2k.
\]
In particular, for every \(m,n\ge 2\), complexities are strict for this word:
\[
R_\alpha(m,n)=m+1< G_\alpha(m,n)=2m< p_\alpha^*(mn)=2mn.
\]

\end{example}
Taking \(r=\min(m,n)\) in \(p_\alpha^*(r)\le G_\alpha(r,r)\le G_\alpha(m,n)\), we obtain the following.
\begin{theorem}
If $\alpha$ is strongly recurrent and $G_\alpha(m,n)<2\min(m,n)$, then $\alpha$ is doubly periodic.
\end{theorem}
\begin{remark}
The inequality $G_\alpha(m,n)<2\min(m,n)$ must be strict; see Example \ref{ex:contex} as a counterexample.
\end{remark}

\section{Simple Toeplitz words}\label{sec:simple}
Let
\[
H_0=\mathbb Z^2 \supsetneq H_1 \supsetneq H_2 \supsetneq \cdots
\]
be a strictly decreasing sequence of finite-index subgroups of $\mathbb Z^2$. 
More precisely, assume that each \(A_r\) is a nonsingular integer matrix with \(|\det(A_r)|\ge 2\); then define $B_0=I$, $B_r=A_1\cdots A_r$ and $H_r=B_r\mathbb Z^2$.

Set $b_0=0$, $D_0=\mathbb Z^2$ and $v_r\in\mathbb Z^2$. For each $r\geq 1$, choose a
vector $b_r=\sum_{j=1}^rB_{j-1}v_j\in\mathbb Z^2$ and a letter $a_r\in\mathcal A$ such that
the coset
\[
D_r=b_r+H_r=\{b_r+h:h\in H_r\}
\]
is contained in $D_{r-1}$.

At stage $r$, assign the same letter $a_r$ to every site in
$D_{r-1}\setminus D_r$. The sites in $D_r$ remain unfilled and form
the \emph{hole coset} at stage $r$. The no-residual-holes condition is
\[
\bigcap_{r\geq 0}D_r=\varnothing.
\]
Thus every site $x\in\mathbb Z^2$ has a unique filling stage
\[
r(x)=\min\{r\geq 1:x\notin D_r\},
\]
and the resulting word is defined by $\alpha(x)=a_{r(x)}$.

A two-dimensional word $\alpha\colon\mathbb Z^2\to\{0,1\}$ is called a \emph{simple Toeplitz word (with no residual holes)} if it admits the above construction and is not doubly periodic.

\begin{theorem}\label{thm:toep}
Let \(\alpha\) be a simple Toeplitz word on \(\mathbb Z^2\). Then
\begin{equation}\label{eq:all}
G_\alpha(m,n)=2mn\quad\text{for every } m,n\ge1\Longleftrightarrow L(\alpha)\cap\mathbb Ze_1=L(\alpha)\cap\mathbb Ze_2=\{0\}.
\end{equation}
\end{theorem}
\begin{lemma}\label{lem:alternate}
	Both letters occur in every tail \(a_{r+1},a_{r+2},\ldots\) and among the
	values of \(\alpha\) on every \(D_r\). After merging consecutive stages
	carrying the same letter, we may assume
	\begin{equation}\label{eq:alternate}
		a_{r+1}=1-a_r\qquad(r\ge1).
	\end{equation}
	Merging preserves the word and both intersections
	\(\bigcap_{r\ge0}D_r\) and \(\bigcap_{r\ge0}H_r\).
\end{lemma}

\begin{proof}
	If the sequence of letters were eventually constant, some \(R\) would
	satisfy \(a_j=a\) for all \(j>R\). From $\cap_{r\geq 0}D_r=\varnothing$, every site of
	\(D_R\) would receive \(a\). Outside \(D_R\), each layer filled by stage
	\(R\) is a union of \(H_R\)-cosets. Thus \(\alpha\) would be constant on
	each \(H_R\)-coset and \(H_R\subset L(\alpha)\), giving double
	periodicity. Consequently both letters occur arbitrarily late. Given
	\(r\), a layer \(D_{j-1}\setminus D_j\) with \(j>r\) is a nonempty subset
	of \(D_r\), so both letters actually occur in \(D_r\).
	
	Each maximal run of identical letters is finite. If the run consists of
	stages \(s+1,\ldots,t\), its union of filled layers is
	\(D_s\setminus D_t\), with a single common letter. Replacing the run by one
	stage therefore preserves \(\alpha\). The retained endpoints tend to
	infinity, so the retained lattices and cosets have the same intersections
	as their original nested chains. Successive runs have different letters,
	which leads to \eqref{eq:alternate}.
\end{proof}

We henceforth use the merged chain and keep its notation. All depths in the
sequel are measured in this chain.

\begin{proposition}\label{prop:per}
	For \(H_\infty=\bigcap_{r\ge0}H_r\), one has
	\begin{equation}\label{eq:per}
	L(\alpha)=H_\infty .
	\end{equation}
\end{proposition}

\begin{proof}
	Every \(v\in H_\infty\) preserves every coset \(D_r\) under translation. It
	therefore preserves every filled layer and its letter. Conversely, take
	\(v\notin H_r\) for some \(r\). The coset \(D_r+v\) is disjoint from \(D_r\).
	Each layer filled by stage \(r\) is a union of \(H_r\)-cosets, so
	\(\alpha\) is constant on \(D_r+v\). By Lemma \ref{lem:alternate}, both
	letters occur on \(D_r\). Translation by \(v\) cannot preserve \(\alpha\).
\end{proof}

For a set \(C\subset\mathbb Z^2\) and a vector \(u\), write
\(C-u=\{e-u:e\in C\}\). For \(a\in\mathcal A\), let \(a^E\) denote the constant
function on \(E\) with value \(a\). Given \(r\ge0\) and a finite nonempty
\(E\subset H_r\), define the \emph{conditional pattern family at level \(r\)}
and its cardinality by
\begin{equation}\label{eq:tail}
	\mathcal F_r(E)=\{(\alpha(b_r+z+e))_{e\in E}:z\in H_r\},
	\quad P_r(E)=|\mathcal F_r(E)|.
\end{equation}
The base point \(x=b_r+z\) ranges over all of \(D_r\), and \(x+E\subset
D_r\). Thus every sampled site is filled after stage \(r\). The choice of
representative \(b_r\) of \(D_r\) has no effect on the family.
In particular, \(\mathcal F_0(E)=\mathcal F_\alpha(E)\) and \(P_0(E)=p_\alpha(E)\).

If \(u\in H_r\), the relabelling \(e\mapsto e-u\) identifies \(\mathcal F_r(E)\)
with \(\mathcal F_r(E-u)\). In fact, a pattern based at \(x\in D_r\) on \(E\)
has the same values as the pattern based at \(x+u\) on \(E-u\), and
\(x\mapsto x+u\) permutes \(D_r\). We call a replacement
\(E\mapsto E-u\) with \(u\in E\) a \emph{normalization}. It places
the origin in the window.

\begin{lemma}\label{lem:before}
	If \(0\le r\le k\) and \(E\subset H_k\) is finite and nonempty, then
	\begin{equation}\label{eq:before}
		\mathcal F_r(E)=\mathcal F_k(E)\cup\{a_j^E:r<j\le k\}.
	\end{equation}
	For \(r=k\), the second set is empty.
\end{lemma}

\begin{proof}
	Partition \(D_r\) into \(D_k\) and the layers \(D_{j-1}\setminus D_j\)
	with \(r<j\le k\). If \(x\) belongs to one of these layers, adding any
	\(e\in E\subset H_k\subset H_j\) preserves both cosets \(D_{j-1}\) and \(D_j\). Hence \(x+e\in D_{j-1}\setminus D_j\). The observed pattern is \(a_j^E\), and this
	pattern is realized since the layer is nonempty. The remaining base points,
	those in \(D_k\), give exactly \(\mathcal F_k(E)\).
\end{proof}

\begin{lemma}\label{lem:split}
	Let \(E\subset H_k\) be finite and nonempty. Partition it into its
	nonempty residue classes modulo \(H_{k+1}\), writing
	\(E=E_1\bigsqcup\cdots\bigsqcup E_s\), where \(\bigsqcup\) denotes disjoint union.
	Choose a representative \(u_i\in E_i\) and put \(F_i=E_i-u_i\subset
	H_{k+1}\). Let \(a=a_{k+1}\) and \(q=[H_k:H_{k+1}]\). Define
	\(J_i\colon\mathcal A^{F_i}\to\mathcal A^E\) by
	\[
	(J_iw)(e)=
	\begin{cases}
		w(e-u_i), & e\in E_i,\\
		a,        & e\in E\setminus E_i .
	\end{cases}
	\]
	Then
	\begin{equation}\label{eq:split}
		\mathcal F_k(E)=\bigcup_{i=1}^s J_i\bigl(\mathcal F_{k+1}(F_i)\bigr)\ \cup B,
		\qquad B=
		\begin{cases}
			\{a^E\},&s<q,\\
			\varnothing,&s=q.
		\end{cases}
	\end{equation}
\end{lemma}

\begin{proof}
	For \(x\in D_k\), a site \(x+e\) is in \(D_{k+1}\) precisely when
	\(e\equiv b_{k+1}-x\pmod{H_{k+1}}\). Thus at most one occupied class of
	the window enters the next hole. All other sites receive the background
	letter \(a\). We call a class that enters the hole the \emph{active class}.
	The base points making \(E_i\) active form exactly the coset
	\[
	b_{k+1}-u_i+H_{k+1}\subset D_k.
	\]
	As \(x\) traverses this coset, \(x+u_i\) traverses all of \(D_{k+1}\).
	Consequently the restrictions to \(E_i\), in the labels of \(F_i\), realize
	every member of \(\mathcal F_{k+1}(F_i)\). Extending by \(a\) on the complement
	gives exactly the indicated embedded family.
	
	As \(x\) ranges over \(D_k\), \(b_{k+1}-x\) traverses all \(q\) residue
	classes of \(H_k/H_{k+1}\). If \(s<q\), an unoccupied class can therefore
	be chosen and the resulting pattern is \(a^E\). If \(s=q\), each
	translation activates one class. This exhausts all possible base points.
\end{proof}
For a finite set \(C\subset\mathbb Z^2\) with at least two points, write
\(C-C=\{x-y:x,y\in C\}\). The condition \(C-C\subset H_k\) means exactly
that all points of \(C\) lie in one \(H_k\)-coset. If the following maximum
is finite, define the \emph{split depth} of \(C\) by
\begin{equation}\label{eq:depth}
	d(C)=\max\{k\ge0:C-C\subset H_k\}.
\end{equation}
Thus \(C\) first separates into distinct classes upon passing from
\(H_{d(C)}\) to \(H_{d(C)+1}\). If \(C-C\subset H_\infty\), we may write
\(d(C)=\infty\), meaning that \(C\) never splits.

The \emph{split tree} of \(E\) has \(E\) as its root. At a node \(C\) of
finite split depth \(k\), its children are its nonempty classes modulo
\(H_{k+1}\). This convention skips intermediate levels at which every point
remains in one class. A singleton is a leaf, and a nonsingleton of infinite
depth is an unsplit terminal node. Translation preserves the difference
set and the residue-class partitions up to relabelling, so it preserves
every split depth.

\begin{lemma}\label{lem:binary}
	Let \(E\subset H_r\) be finite and nonempty. Suppose its split tree
	satisfies the following three conditions:
	\begin{enumerate}
		\renewcommand{\labelenumi}{\textup{(\roman{enumi})}}
		\item Every nonsingleton node splits at finite depth into exactly two
		nonempty children.
		\item If the root is a nonsingleton, its split depth is at least \(r+2\).
		\item For each nonsingleton child, its split depth is at least the split
		depth of its parent plus three.
	\end{enumerate}
	Then
	\begin{equation}\label{eq:exact}
		P_r(E)=2|E|,\qquad \{0^E,1^E\}\subset\mathcal F_r(E).
	\end{equation}
\end{lemma}

\begin{proof}
	We induct on \(|E|\), allowing the initial level \(r\) to vary. For
	\(E=\{e\}\), the site \(x+e\) ranges over all of \(D_r\) as \(x\) ranges
	over \(D_r\). Both letters occur there by Lemma~\ref{lem:alternate}.
	This proves both conclusions at a leaf.
	
	Let \(k\) be the split depth of a nonsingleton root. Choose \(u\in E\)
	and replace \(E\) by \(E-u\). This preserves its conditional family up
	to relabelling at level \(r\), and \(E-u\subset H_k\) by
	\eqref{eq:depth}. Denote the normalized window again by \(E\). Its first
	split is \(E=E_1\bigsqcup E_2\). Put \(N_i=|E_i|\) and take
	\(F_i=E_i-u_i\subset H_{k+1}\) as in Lemma~\ref{lem:split}.
	For each child, the conditional family now starts at level \(k+1\).
	If the child is a nonsingleton, its first split has depth at least
	\(k+3=(k+1)+2\). All deeper parent--child gaps are inherited. The induction
	hypothesis therefore gives, for \(i=1,2\),
	\[
	|\mathcal F_{k+1}(F_i)|=2N_i,\qquad
	\{0^{F_i},1^{F_i}\}\subset\mathcal F_{k+1}(F_i).
	\]
	
	Set \(a=a_{k+1}\). The two embedded child families both contain \(a^E\), obtained by choosing \(a^{F_i}\) on the active child. Their intersection consists of
	this one pattern, since a nonbackground site cannot lie in both disjoint
	children. If some translation activates neither child, it adds the same
	constant pattern. Hence \eqref{eq:split} gives
	\begin{equation}\label{eq:minusone}
		P_k(E)=2N_1+2N_2-1=2|E|-1 .
	\end{equation}
	The other constant pattern \((1-a)^E\) is absent from \(\mathcal F_k(E)\) because each translation leaves at least one nonempty child entirely in the
	background. Since \(k\ge r+2\), the alternating list
	\(a_{r+1},\ldots,a_k\) contains both letters. By \eqref{eq:before}, the
	earlier base points contribute precisely the missing constant pattern in
	addition to those already counted. Thus \(P_r(E)=2|E|\) and both
	constant patterns are realized.
\end{proof}
\begin{proof}[Proof of Theorem \ref{thm:toep}]
(Sufficiency) For \(v\in\mathbb Z^2\setminus H_\infty\) define its \emph{depth}
\(\ell(v)=\max\{r\ge0:v\in H_r\}\). Set \(\ell(v)=\infty\) for
\(v\in H_\infty\), so \(\ell(0)=\infty\). We order \(\infty\) above
all finite depths. As the \(H_r\) are subgroups, \(\ell(-v)=\ell(v)\).
For vectors of unequal depths we have
\begin{equation}\label{eq:nocancel}
	\ell(v)<\ell(w)\quad\Longrightarrow\quad
	\ell(v+w)=\ell(v-w)=\ell(v).
\end{equation}
Indeed, put \(h=\ell(v)\). Both vectors lie in \(H_h\), while modulo
\(H_{h+1}\) the vector \(v\) is nonzero and \(w\) is zero. The two sums
\(v\pm w\) thus have nonzero residue at precisely this first distinguishing
level.

For \(L\ge0\) and \(i\in\{1,2\}\), define
\begin{equation}\label{eq:qi}
	q_i(L)=\min\{q\ge1:q e_i\in H_L\}.
\end{equation}
This is the order of the coset \(e_i+H_L\) in the finite quotient group
\(\mathbb Z^2/H_L\), so the minimum exists. Suppose that neither coordinate direction is periodic, where \(e_1=(1,0)\) and \(e_2=(0,1)\). By \eqref{eq:per}, \(q_i(L)e_i\notin H_\infty\), so its depth is finite and at least \(L\).

Starting with \(L=2\), choose \(s_1,\ldots,s_{m-1}\) on the first axis, then \(t_1,\ldots,t_{n-1}\) on the second. On axis \(i\), take the next coordinate to be \(q_i(L)\) and reset \(L=\ell(q_i(L)e_i)+3\). Skip an empty list. Thus the chosen depths increase by at least three.

Let
\begin{equation}\label{eq:window}
	S=\{0,s_1,\ldots,s_{m-1}\},\qquad
	T=\{0,t_1,\ldots,t_{n-1}\},\qquad E=S\times T.
\end{equation}
By \eqref{eq:nocancel}, differences between distinct horizontal coordinates have shallower depths and cannot be cancelled by vertical differences. When $m>1$, the first split therefore separates the entire group with horizontal coordinate $s_1$ from the rest. Within the remaining group, the coordinates $s_2,s_3,\ldots$ are separated successively. Inside a group whose horizontal coordinate is already fixed, the coordinates $t_1,t_2,\ldots$ split successively. Every nonsingleton node has exactly two children, the leaves are
singletons, and the splitting depths meet every hypothesis of Lemma \ref{lem:binary}. Hence
\[p_\alpha(E)=2|E|=2mn.\]
Together with $G_\alpha(m,n)\leq p_\alpha^*(mn)=2mn$~\cite{qrwx}, this proves $G_\alpha(m,n)=2mn$.

(Necessity) If $(0,q)$ is a period, where $q\geq 1$, every pattern on $S \times T$ is determined by representatives of the distinct vertical coordinates modulo $q$. These representatives give at most $m\min(n,q)$ sites, so
\[
G_\alpha(m,n)\leq p_\alpha^*(m\min(n,q))=2m\min(n,q).
\]
For $n>q$ this is strictly less than $2mn$. The argument for a horizontal period is identical.
\end{proof}
We straightforwardly have the following corollaries.
\begin{corollary}
Let \(\alpha\) be a simple Toeplitz word on \(\mathbb Z^2\).  
\begin{equation}\label{eq:some}
	\begin{aligned}
			&G_\alpha(m,n)<2mn\quad\text{for some }m,n\ge1\\
			&\Longleftrightarrow
			\alpha\text{ has a horizontal or a vertical period}.
	\end{aligned}
\end{equation}
\end{corollary}
\begin{corollary}
Let \(\alpha\) be a simple Toeplitz word. If $(0, q)$ is a period with \(q\ge1\), then
\[G_\alpha(m, n) \leq  2m\min(n, q),\quad G_\alpha(m, n) = G_\alpha(m, q)\ (n \geq q). \]
If $(p, 0)$ is a period with \(p\ge1\), then
\[G_\alpha(m, n) \leq  2n\min(m, p),\quad G_\alpha(m, n) = G_\alpha(p,n)\ (m \geq p). \]
\end{corollary}
\begin{example}
On \(\mathbb Z^2\), the word of Kamae and Xue \cite{kx} is given by
\[
\alpha(i,j)=\begin{cases}
	1, &\text{if} \deg_2(i)=\deg_2(j),\cr
	0, &\text{otherwise},
\end{cases}
\]
where, for \(i\ne0\), $\deg_2(i)$ is the largest integer $r$ such that $2^r\mid i$. Set $\deg_2(0)=\infty$.

Let $a_{2k-1}=0$, $a_{2k}=1$, $v_k=0$ and 
\[A_k=A=\begin{bmatrix}1&1\\
1&-1\end{bmatrix}\]
for all $k\in\mathbb Z_+$. Here \(D_r=H_r=A^r\mathbb Z^2\) and \(A^2=2I\). Thus \(H_{2k}=2^k\mathbb Z^2\), and \(H_{2k+1}\) consists of \(2^k(u,v)\) with \(u\equiv v\pmod 2\). The word takes the stated alternating letters on the successive layers. The origin remains in every \(H_r\), so we use the split count directly. Indeed, the proof of Lemma~\ref{lem:binary} uses only the layer decompositions \eqref{eq:before} and \eqref{eq:split} and the presence of both letters in every hole coset. The decompositions hold here by the same finite coset partitions, and every \(H_r\) contains layers of both letters.

Take \(S=\{0\}\cup\{4^j:1\le j<m\}\) and \(T=\{0\}\cup\{4^{m-1+j}:1\le j<n\}\). Since \(\ell(4^je_i)=4j\) for either axis, the split argument for \eqref{eq:window} yields \(p_\alpha(S\times T)=2mn\).
For the upper bound, every finite pattern occurs in the positive quadrant: translating by \((2^M,2^M)\) preserves all its letters for large \(M\), since each nonzero coordinate keeps its valuation and each zero coordinate acquires valuation \(M\). The bound \(p_\alpha^*(mn)\le 2mn\) of \cite{kx} therefore applies to the Kamae--Xue word on $\mathbb Z^2$. Hence
\[G_\alpha(m,n)=2mn,\quad m,n\geq 1.\]
\end{example}
\section{Irrational rotation words}\label{sec:irr}
Let $\mathbb T=\mathbb R/\mathbb Z$, $\lambda\in(0,1)$, $J=[0,\lambda)\ (\bmod\ 1)$, and
\[
\alpha(i,j)=\bm{1}_J(t+ia+jb\ (\bmod 1)).
\]
Unless otherwise specified, we assume throughout this section that at least one of $a,b$ is irrational. The translation phases $t+a\mathbb Z+b\mathbb Z$ are then dense in $\mathbb T$. For rectangular factors of rotation words, see Berth\'e and Vuillon \cite{bv}.

\begin{lemma}\label{lem:upper}
For a window $E$ of size $N$, we have
\[
p_\alpha(E)\leq 2N.
\]
\end{lemma}
\begin{proof}
Each of the \(N\) sites contributes two endpoints of a translated interval \(J\). The pattern is constant on each of the resulting at most \(2N\) arcs, and the half-open convention gives each endpoint the pattern immediately to its right.
\end{proof}
\begin{lemma}\label{lem:2N}
If the phases of the $N$-sites in the window  have distinct real representatives $0 \leq \delta_1 <\cdots  < \delta_N < \varepsilon$, with $\varepsilon < \min(\lambda, 1 -\lambda)$, the window has exactly $2N$ patterns.
\end{lemma}
\begin{proof}
Index the window sites by their phases
	\(\delta_1<\cdots<\delta_N\). At translation phase \(x\), coordinate
	\(k\) equals \(1\) precisely when
	\[
	x\in[A_k,B_k)\pmod 1,\quad
	A_k=-\delta_k,\quad B_k=\lambda-\delta_k.
	\]
	Put \(D=\delta_N-\delta_1\). Since \(D<\min\{\lambda,1-\lambda\}\),
	the gaps \(B_N-A_1=\lambda-D\) and \(A_N+1-B_1=1-\lambda-D\)
	are both positive. Thus, over one turn, the endpoints occur in the order
	\[
	A_N<\cdots<A_1<B_N<\cdots<B_1<A_N+1.
	\]
	
	Starting with empty support just before \(A_N\), the coordinates turn
	on in the order \(N,\ldots,1\), then turn off in the same order.
	The supports are therefore the empty set, the full set, the \(N-1\)
	nonempty proper suffixes of \(\{1,\ldots,N\}\), and its \(N-1\)
	nonempty proper prefixes. Within each family the sizes are distinct,
	and every proper suffix contains \(N\) whereas every proper prefix
	omits it. This gives exactly \(2N\) distinct binary patterns.
	
	Each pattern persists on an arc of positive length and is realized
	because the translation phases are dense. At each endpoint the
	half-open convention gives the pattern immediately to its right.
	Hence \(p_\alpha(E)=2N\).
\end{proof}
\begin{theorem}
If $a,b$ are irrational, then $G_\alpha(m,n)=2mn$ for $m,n\geq 1$.
\end{theorem}
\begin{proof}
To choose the grid, assume that both \(a\) and \(b\) are irrational.
Fix \(m,n\ge1\) and put \(\varepsilon=\tfrac12\min\{\lambda,1-\lambda\}\).
Writing \(\{x\}\) for the fractional part of \(x\), we first choose \(v\ge1\)
and then \(u\ge1\) such that
\[
0<\beta:=\{vb\}<\frac{\varepsilon}{n},
\quad
0<\gamma:=\{ua\}<\frac{\beta}{m}.
\]
Both choices are possible by density of the positive multiples of each
irrational number modulo one. Take the grid
\[
S=\{0,u,\ldots,(m-1)u\},\qquad
T=\{0,v,\ldots,(n-1)v\}.
\]

The observation phases in row \(j\), for \(0\le j<n\), have representatives
\[
j\beta, j\beta+\gamma,\ldots, j\beta+(m-1)\gamma.
\]
They increase strictly within the row. Since \((m-1)\gamma<\beta\),
the last phase in row \(j\) lies before \((j+1)\beta\), the first phase
of the next row. Thus all \(mn\) phases are distinct, and the largest is
\[
(n-1)\beta+(m-1)\gamma<n\beta<\varepsilon.
\]
Lemma \ref{lem:2N} now gives \(p_\alpha(S\times T)=2mn\). Lemma \ref{lem:upper} shows the
matching upper bound for every \(m\times n\) grid.  So
\(G_\alpha(m,n)=2mn\).
\end{proof}
\begin{remark}
Rational independence of $1, a, b$ is not needed.
\end{remark}
\section{Words with full complexity}\label{sec:full}
We now construct two two-dimensional binary words with full grid complexity.
\subsection{A word with 2-adic construction}
\begin{proposition}
	For each positive integer \(x\), let \(v_2(x)\) denote its
	\(2\)-adic valuation, defined by
	\[
	v_2(x)=\max\{k\in\mathbb Z_{\ge 0}:2^k\mid x\}.
	\]
	Thus \(x=2^{v_2(x)}o(x)\), where \(o(x)\) is the unique odd
	positive integer in this factorization.
	Let
	\[
	d_k(r)=\left\lfloor \frac{r}{2^k}\right\rfloor\pmod 2
	\]
	denote the \(k\)-th binary digit of \(r\), starting at \(k=0\).
	Define a binary word on \(\mathbb Z^2\) by
	\[
	\alpha(x,y)=
	\begin{cases}
		d_y(o(x)), & x>0,\ y\ge0,\\
		0,         & \text{otherwise}.
	\end{cases}
	\]
	Then, for all \(m,n\ge1\),
	\[
	G_\alpha(m,n)=p_\alpha^*(mn)=2^{mn}.
	\]
	In fact, the grid
	\[
	S\times T
	=\{2^{in}:0\le i<m\}\times\{0,1,\ldots,n-1\}
	\]
	realizes all \(2^{mn}\) binary patterns.
\end{proposition}

\begin{proof}
	Fix \(m,n\ge1\), put \(L=mn\), and let
	\(B=(b_{ij})_{0\le i<m,\,0\le j<n}\) be any prescribed binary matrix.
	Encode its rows in the integer
	\[
	M=\sum_{i=0}^{m-1}\sum_{j=0}^{n-1} b_{ij}\,2^{in+j},
	\quad\text{so that}\quad d_{in+j}(M)=b_{ij}.
	\]
	Translate the grid by \((2^L M,L)\). At its \(i\)-th horizontal
	coordinate we have
	\[
	2^L M+2^{in}=2^{in}\bigl(1+2^{L-in}M\bigr).
	\]
	Since \(L-in\ge n\ge1\), the factor in parentheses is odd.
	Its binary expansion is therefore obtained by shifting the digits
	of \(M\) to the left by \(L-in\) places and setting the lowest bit
	to \(1\). In particular, its digit at position \(L+j\) is the
	digit of \(M\) at position \(in+j\). Hence
	\[
	\alpha(2^L M+2^{in},L+j)
	=d_{L+j}\bigl(1+2^{L-in}M\bigr)
	=d_{in+j}(M)
	=b_{ij}.
	\]
	Thus every \(B\) occurs on a translate of the same grid \(S\times T\).
	Since any binary window of \(mn\) sites admits at most \(2^{mn}\)
	patterns, both equalities follow.
\end{proof}
\subsection{Pascal word}
\begin{definition}
	The binary Pascal word \(P:\mathbb Z^2\to\{0,1\}\) is defined by
	\[
	P(x,y)=
	\begin{cases}
		\displaystyle\binom{x+y}{x}\bmod 2, & x,y\ge 0,\\[4pt]
		0, & \text{otherwise}.
	\end{cases}
	\]
\end{definition}

\begin{proposition}
	For every \(m,n\ge 1\), the Pascal word satisfies
	\[
	G_P(m,n)=p_P^*(mn)=2^{mn}.
	\]
	More precisely, there is an \(m\times n\) grid on which horizontal
	translations alone realize every binary \(m\times n\) pattern.
\end{proposition}

\begin{proof}
	For \(x,y\ge0\), recall that \(P(x,y)=1\) precisely when \(x\) and \(y\)
	have no common \(1\) in their binary expansions. Indeed, if
	\(x+y=\sum_{k\in K}2^k\), then over \(\mathbb F_2\),
	\[
	(X+Y)^{x+y}=\prod_{k\in K}(X^{2^k}+Y^{2^k}).
	\]
	The coefficient of \(X^xY^y\) is \(1\) exactly when the binary
	digits indexed by \(K\) can be partitioned between \(x\) and \(y\),
	which gives the stated criterion.
	
	Fix \(m,n\ge 1\). Assign to each cell \((r,c)\) the base-four
	position \(rn+c\), where \(0\le r<m\) and \(0\le c<n\). Take
	\[
	S=\{s_i:0\le i<m\},\qquad s_i=\sum_{c=0}^{n-1}4^{in+c},
	\quad
	T=\{t_j:0\le j<n\},\qquad t_j=2\sum_{r=0}^{m-1}4^{rn+j}.
	\]
	Thus \(s_i\) has digit \(1\) in the positions assigned to row \(i\),
	and \(t_j\) has digit \(2\) in those assigned to column \(j\).
	All their other digits are zero.
	
	Given a target matrix \(B=(b_{rc})\in\{0,1\}^{m\times n}\), set
	\[
	u_B=\sum_{r=0}^{m-1}\sum_{c=0}^{n-1}(1-b_{rc})\,4^{rn+c}.
	\]
	The base-four digits of \(u_B\) are \(0\) or \(1\), so adding \(s_i\)
	causes no base-four carry. The resulting digit equals \(2\) exactly
	at the cells \((i,c)\) for which \(b_{ic}=0\). Since \(2=(10)_2\),
	these are precisely the positions where the higher bit of a binary
	digit pair in \(u_B+s_i\) is \(1\). The integer \(t_j\) tests that
	higher bit at every cell of column \(j\). Hence \(u_B+s_i\) and
	\(t_j\) have a common binary \(1\) exactly when \(b_{ij}=0\).
	By the parity criterion,
	\[
	P(u_B+s_i,t_j)=b_{ij}
	\qquad(0\le i<m,\ 0\le j<n).
	\]
	
	Thus the horizontal translate \((u_B,0)+(S\times T)\) realizes \(B\).
	Since the grid is fixed and \(B\) is arbitrary,
	\[
	2^{mn}\le G_P(m,n)\le p_P^*(mn)\le 2^{mn},
	\]
	as claimed.
\end{proof}

\section*{Acknowledgments}
This work was partially supported by the National Natural Science Foundation of China (Nos. 12301107, 62450131, and 12271038) and the Shandong Province Natural Science Foundation, China (No. ZR202209010046).

\section*{Author Contributions}
Cheng Zeng: Methodology; Formal analysis; Conceptualization; Funding acquisition; Writing -- original draft.
Qingcheng Zeng: Formal analysis;  Writing -- review \& editing; Validation. 
Yumei Xue: Writing -- review \& editing; Supervision; Funding acquisition. 
All authors have read and agreed to the published version of the manuscript.

\section*{Conflicts of Interest}
The authors declare no conflicts of interest.

\section*{Use of generative AI} 
During preparation of this manuscript, the authors used ChatGPT to assist with language revision and literature searches. Responsibility for the mathematical arguments, citations, and final text rests with the authors.

\end{document}